\documentclass{article}
\title{\large \bf On central automorphisms of finite  $p$-groups }
\author{\small \bf Deepak Gumber \\
\small \em School of Mathematics and Computer Applications \\
\small \em Thapar University, Patiala, India\\
\small \em Email: dkgumber@yahoo.com\\
\\
\small \bf Mahak Sharma\\
\small \em Department of Applied Sciences and Humanities\\
\small \em Baddi University of Emerging Sciences and Technology, Baddi, India \\
\small \em Email: mathes4282@gmail.com}

\date{}
\newtheorem{thm}{Theorem}[section]
\newtheorem{lm}[thm]{Lemma}

\usepackage{amssymb}

\begin{document}
\maketitle
\begin{abstract}
We characterize all finite $p$-groups $G$ of order $p^n(n\leq 6)$, where $p$ is a prime for $n\leq 5$ and an odd prime for $n=6$, such that the center of the inner automorphism group of $G$ is equal to the group of central automorphisms of $G$.
\end{abstract}

\vspace{2ex}


\vspace{2ex}


\section{\large Introduction}
Let $G$ be a finite $p$-group. By $Z(G),\gamma_2(G),
\mbox{Aut}(G)$ and $\mbox{Inn}(G)$, we denote the
center, the commutator subgroup, the group of all
 automorphisms and the group of all
inner automorphisms of $G$ respectively. An automorphism
$\alpha $ of $G$ is called a central automorphism if
 $\alpha$ commutes with every inner automorphism, or equivalently, if $x^{-1}
\alpha(x)$ lies in the centre of $G$
 for all $x$ in $G$. The central automorphisms fix the
 commutator subgroup $\gamma_{2}(G)$ of $G$
pointwise and form a normal subgroup $\mbox{Aut}_z(G)$ of $\mbox{Aut}(G)$. The group $\mathrm{Aut}_z(G)$ always lies between $Z(\mathrm{Inn}(G))$ and $\mathrm{Aut}(G)$. Non-abelian $p$-groups in which all automorphisms are central have been well studied (see [6] for references). Curran and McCaughan [6] have characterized those finite $p$-groups for which $\mathrm{Aut}_z(G)=\mathrm{Inn}(G)$. They in fact proved that
$\mbox{Aut}_z(G)=\mbox{Inn}(G)$ if and only if $\gamma_2(G)=Z(G)$ and $Z(G)$ is cyclic. Curran [5] considered the case when $\mathrm{Aut}_z(G)$ is minimum possible, that is, when $\mathrm{Aut}_z(G)=Z(\mathrm{Inn}(G))$. He proved that for $\mbox{Aut}_z(G)$ to be equal to $Z(\mbox{Inn}(G))$, $Z(G)$ must be contained in $\gamma_2(G)$ and $Z(\mbox{Inn}(G))$ must not be cyclic. These conditions are necessary but not sufficient because there are groups $G$ for which $\mbox{Aut}_z(G)\neq Z(\mbox{Inn}(G))$ even if $Z(G)$ is contained in $\gamma_2(G)$ or $Z(\mbox{Inn}(G))$ is not cyclic [see section 4]. The present note is the result of an effort to find the sufficient conditions. Observe that if nilpotency class of $G$ is 2, then $Z(\mbox{Inn}(G))=\mbox{Inn}(G)$ and thus $\mbox{Aut}_z(G)=Z(\mbox{Inn}(G))$ if and only if $\gamma_2(G)=Z(G)$ and $Z(G)$ is cyclic. Also observe that if $G$ is of maximal class, then $|Z(\mbox{Inn}(G))|=p$ and by [5] $\mbox{Aut}_z(G)>Z(\mbox{Inn}(G))$. Therefore, to characterize all finite $p$-groups $G$ for which $\mbox{Aut}_z(G)=Z(\mbox{Inn}(G))$, we can assume that nilpotency class of $G$ is bigger than 2 and $G$ is not of maximal class. Of course in such a case $|G|\geq p^5$. In section 3, we  prove two theorems which characterize all such groups of order $p^5$, where $p$ is any prime, and of order $p^6$, where $p$ is an odd prime. In section 4, we find all such groups from the list of groups of order $p^n$, where $p$ is an odd prime and $4\leq n\leq 6$, ordered into isoclinism families by James [7].

\section{\large Preliminaries} By $\mathrm{Hom}(G,A)$, we denote the group of all homomorphisms of $G$ into an abelian group $A$ and by $C_{p^n}$, we denote the cyclic group of order $p^n$. The nilpotency class of $G$ is denoted as $cl(G)$ and by $d(G)$ we denote the rank of $G$.  For a $p$-group $G$, let $\Omega_1(G)=\;<x\in G|x^p=1>$ and
$\mho_1(G)=\;<x^p|x\in G>$.  A non-abelian group $G$ that has no non-trivial abelian direct factor is said to be purely non-abelian. Observe that a group $G$ is purely non-abelian if its center $Z(G)$ is contained in the frattini subgroup $\Phi(G)$. A $p$-group $G$ is said to be regular if for any two elements $x,y$ in $G$, there is an element $z$ in the commutator subgroup $\gamma_2(H)$ of the subgroup $H=\;<x,y>$ of $G$, such that $x^py^p=(xy)^pz^p$. In the following lemma we list two important properties of the groups $G$ for which 
$\mathrm{Aut}_z(G)=Z(\mathrm{Inn}(G))$.

\begin{lm}{\em [5, Corollaries 3.7, 3.8]} Let $G$ be a finite non-abelian $p$-group such that $\mathrm{Aut}_z(G)=Z(\mathrm{Inn}(G))$. Then $Z(G)\le \gamma_2(G)$ and $Z(\mathrm{Inn}(G))$ is not cyclic.
\end{lm}

The following well known results are also of our interest and will be referred to without any citation. 

\begin{thm}{\em [1]}
If $G$ is a purely non-abelian finite group,
then $|\mathrm{Aut}_z(G)|=|\mathrm{Hom}(G/\gamma_2(G),Z(G))|$.
\end{thm}

\begin{lm}
Let A,B and C be finite abelian groups. Then
$$(i)\;\; \mathrm{Hom}(A \times B, C) \cong \mathrm{Hom}(A,C) \times \mathrm{Hom}(B,C).$$
$$(ii)\;\; \mathrm{Hom}(A,B \times C) \cong \mathrm{Hom}(A,B) \times \mathrm{Hom}(A,C).$$
\end{lm}
\begin{lm}
 Let $C_m$ and $C_n$ be
 two cyclic groups of order $m$ and $n$ respectively. Then
$ \mathrm{Hom}(C_m,C_n) \cong C_{d} $ ,
 where d is the greatest common divisor of  m and n.
\end{lm}

\begin{lm}
If $G$ is a regular $p$-group, then for all $x,y\in G$ and $i,j\geq 0$,
$$[x^{p^i},y^{p^j}]=1\;\mbox{if and only if}
\;\;[x,y]^{p^{i+j}}=1.$$
\end{lm}

\section{\large Characterization}
In this section we characterize all groups of order $p^5$, where $p$ is any prime, and all groups of order $p^6$, where $p$ is an odd prime, such that $\mathrm{Aut}_z(G)=Z(\mathrm{Inn}(G))$. The characteristics of such groups turn up same in both the cases.
We begin with the following simple lemma.
\begin{lm} 
Let $G$ be a finite $p$-group of order $p^4$ and rank $2$. If $|\gamma_2(G)|=p$, then $Z(G)=\Phi(G)$ and $|Z(G)|=p^2$.
\end{lm}
{\bf Proof.} Observe that $|\Phi(G)|=p^2$ and $cl(G)=2$. Therefore $\exp(G/Z(G))=\exp(\gamma_2(G))=p$ and hence $Z(G)=\Phi(G)$.

\begin{thm}
Let $G$ be a finite $p$-group such that $|G|=p^5$ and $cl(G)=3$. Then
$\mathrm{Aut}_z(G)=Z(\mathrm{Inn}(G))$ if and only if $d(G)=2$ and $|Z(G)|=p.$
\end{thm}
{\bf Proof.} First suppose that $d(G)=2$ and $|Z(G)|=p$. Then $G$ is purely non-abelian and thus $|\mathrm{Aut}_z(G)|=|\mathrm{Hom}(G/\gamma_2(G),Z(G))|
=p^2.$
Since $Z(\mathrm{Inn}(G))\leq \mathrm{Aut}_z(G)$, $|Z(\mathrm{Inn}(G))|=p$ or $p^2$. If $|Z(\mathrm{Inn}(G))|=p$, then $G/Z(G)$ is a 2-generated non abelian group of order $p^4$ having nilpotency class 2. Thus $|Z(G/Z(G))|=|Z(\mathrm{Inn}(G))|=p=\gamma_2(G/Z(G))$, which is a contradiction to Lemma 3.1. Therefore $|Z(\mathrm{Inn}(G))|=p^2$ and hence $\mathrm{Aut}_z(G)=Z(\mathrm{Inn}(G))$.

Conversely suppose that $\mbox{Aut}_z(G)=Z(\mbox{Inn}(G))$. Then  $Z(G)<\gamma_2(G)$ by Lemma 2.1 and the fact that $cl(G)=3$. Therefore  $G$ is purely non-abelian and
$\mathrm{Inn}(G)$ is a non-abelian group of order $p^4$, because if $|\mathrm{Inn}(G)|=p^3$ , then $|Z(\mbox{Inn}(G))|=p$ which is a contradiction to Lemma 2.1. Thus  $|Z(G)|=p$ and    
$|Z(\mbox{Inn}(G))|=p^2$. Then 
$|\mathrm{Aut}_z(G)|=|\mathrm{Hom}(G/\gamma_2(G),Z(G))|=
|Z(\mbox{Inn}(G))|=p^2$ implies that $d(G)=2$.

\begin{thm}
Let $G$ be a finite $p$-group, $p$ an odd prime, such that $|G|=p^6$ and $cl(G)=3$ or $4$. Then
$\mathrm{Aut}_z(G)=Z(\mathrm{Inn}(G))$ if and only if $d(G)=2$ and $|Z(G)|=p.$
\end{thm}
{\bf Proof.} If $d(G)=2$ and $|Z(G)|=p$, then as in above theorem,  $|\mbox{Aut}_z(G)|=p^2$. We proceed to prove that $|Z(\mbox{Inn}(G))|=p^2$. On the contrary suppose that $|Z(\mbox{Inn}(G))|=p$. First suppose that $cl(G)=3$. Let $H=G/Z(G)$. Since the nilpotency class of $G$ is 3, the nilpotency class of $H$ is 2. Thus $\gamma_2(H)=Z(H)=Z(\mathrm{Inn}(G))$ is of order $p$. This implies that the exponent of $H/Z(H)$ is $p$ and therefore $H$ is an extra-special $p$-group, which is a contradiction to the fact that $d(H)=2$. Next suppose that  $cl(G)=4$. Let $H=G/Z(G)$ and let $K=H/Z(H)$. Then $|K|=p^4$, $cl(K)=2$ and $\exp(K/Z(K))=\exp(\gamma_2(K))$. Now $|Z(K)|\neq p$ because if $|Z(K)|=p$, then $d(K)=3$. Thus $|Z(K)|=p^2$ and hence $|\gamma_2(K)|=p$ [8, Theorem 2.1]. It then follows that $|\gamma_2(H)|=p^2,\;|\gamma_2(G)|=p^3$ and $|\gamma_3(G)|=p^2$. Let $\{x,y\}$ be a minimal generating set for $H$. If $\gamma_2(H)$ is elementary abelian, then
$$[x^p,y]=[x,y]^{x^{p-1}}[x,y]^{x^{p-2}}\ldots [x,y]^x[x,y]=[x,y]^p[x,y,x]^{p(p-1)/2}=1.$$
Thus $x^p$ and similarly $y^p$ is in $Z(H)<\gamma_2(H)$. Therefore $\Phi(H)=\gamma_2(H)$, a contradiction to the fact that $d(H)=2$. If $\gamma_2(H)$ is cyclic, then  $|\Phi(\gamma_2(H))|=p$. But $\Phi(\gamma_2(H))=\Phi(\gamma_2(G)/Z(G))=
\Phi(\gamma_2(G))/Z(G)$, because $Z(G)$ being of order $p$ is contained in $\Phi(\gamma_2(G))$. This implies that $|\Phi(\gamma_2(G))|=p^2$. It then follows that $\gamma_2(G)$ is cyclic and therefore $G$ is regular. Suppose
$G/\gamma_2(G)\approx C_{p^2}\times C_p\approx \;<x\gamma_2(G)>\times <y\gamma_2(G)>,$ then
$\{x,y\}$ is a minimal generating set for $G$, $\gamma_2(G)=\;<[x,y]>$ and $x^{p^2},y^p\in \gamma_2(G)$ but $x^p\notin \gamma_2(G)$. Since $G$ is of class 4, $\gamma_3(G)\leq Z_2(G)$. But $|Z_2(G)/Z(G)|=|Z(\mathrm{Inn}(G))|=p=|Z(G)|$ implies that $|Z_2(G)|=p^2=|\gamma_3(G)|$ and thus $\gamma_3(G)=Z_2(G)$. If $y^p\in Z_2(G)=\gamma_3(G)$, then $[y^p,x]\in Z(G)$ and since $G$ is regular, $1=[y^p,x]^p=[y^{p^2},x]=[y,x]^{p^2}$, a contradiction to the fact that $|\gamma_2(G)|=p^3$. Therefore $y^p\in \gamma_2(G)-\gamma_3(G)$, $\gamma_2(G)=\;<y^p>$ and $|y|=p^4$.
Thus $x^{p^2}\in\; <y>$ and therefore $[x^{p^2},y]=[x,y]^{p^2}=1$, which is again a contradiction to the fact that $|\gamma_2(G)|=p^3$. Hence $\mathrm{Aut}_z(G)=Z(\mathrm{Inn}(G))$ is of order $p^2$.

Conversely suppose that $Z(\mbox{Inn}(G))=\mathrm{Aut}_z(G)$. Then $Z(G)<\gamma_2(G)$ by Lemma 2.1 and the fact that $cl(G)>2$. Thus $G$ is purely non-abelian and $G/Z(G)$ is a non-abelian group of order at most $p^5$. Therefore $p^2\leq|Z(\mbox{Inn}(G))|\leq p^3$, because $Z(\mathrm{Inn}(G))$ cannot be cyclic by Lemma 2.1. The possibility of $G$ being rank 4 is immediately ruled out because if $d(G)=4$, then
$$|\mathrm{Aut}_z(G)|=|\mathrm{Hom}({C_{p}\times C_{p}\times C_{p}\times C_{p},C_{p}})|=p^{4}\neq |Z(\mbox{Inn}(G))|.$$
If $d(G)=3$, then $|\mathrm{Aut}_z(G)|\geq p^3$ and therefore by assumption, 
$|\mathrm{Aut}_z(G)|=|Z(\mbox{Inn}(G))|=p^3.$
Also $|Z(G)|= p$, for if  $|Z(G)|=p^2$, then $|Z(\mathrm{Inn}(G))|= p^2$. Let $H=G/Z(G)$, then $|H/Z(H)|=p^2$ and hence $|\gamma_2(H)|=p$. It then follows that  $|\gamma_2(G)|=p^2,\;|\gamma_3(G)|=p$ and therefore $cl(G)=3$. Suppose that $$G/\gamma_2(G)\approx C_p\times C_p\times C_{p^2}\approx \;<x\gamma_2(G)>\times <y\gamma_2(G)>\times <z\gamma_2(G)>.$$ Then $G=\;<x,y,z>$ and $x^p,y^p,z^{p^2}\in \gamma_2(G)$ but $z^p\notin \gamma_2(G).$  Since $|G/Z_2(G)|=p^2$, one of $x,y$ and $z$ lies in $Z_2(G)$. If $z\in Z_2(G)$, then $z^p\in Z(G)\leq \gamma_2(G)$, which is a contradiction. We can therefore, without any loss of generality, assume that $x\in Z_2(G).$ Then $[x,y],[x,z]\in Z(G)$ and $[z^p,x]=1$. If $\gamma_2(G)$ is elementary abelian, then 
$[z^p,y]=[z,y]^p[z,y,z]^{p(p-1)/2}=1.$
Thus $z^p\in Z(G)$, again a contradiction. Now suppose that $\gamma_2(G)$ is cyclic. Let $M=\;<y,z,\Phi(G)>,$ then $M$ is a maximal subgroup of $G$ and $G=MZ_2(G)$. It follows from [3, Theorem 1.3] that $\gamma_i(G)=\gamma_i(M)$ for all $i\ge 2$. Since $Z(G)\leq M$, $Z(G)\le Z(M)$ and $C_G(M)=Z(M)$. We prove that $Z(G)=Z(M)$. Since $|M|=p^5$, order of $Z(M)$ lies between $p$ and $p^3$. If $|Z(M)|=p^2$ or $p^3$, then since $\gamma_2(M)=\gamma_2(G)$ is cyclic, it follows from [2, Prop. 21.20] that 
$$p^3\ge |M/Z(M)|\ge |\gamma_2(M)|^2=p^4,$$
a contradiction and thus $|Z(M)|=p=|Z(G)|$. This proves that $Z(M)=Z(G)$. Since $\gamma_2(G)$ is cyclic, $\Omega_1(\gamma_2(G))=Z(G)$ and thus $C_G(C_G(M))=M$ [4, Theorem C]. Now from the facts that $C_G(M)=Z(M)$, $Z(M)=Z(G)$ and $C_G(C_G(M))=M$, it follows that $C_G(Z(G))=M$, a final contradiction to $d(G)=3$. 
Hence $d(G)$=2 and therefore, because of assumption, $|Z(G)|=p$ or $p^2$. We prove that  $|Z(G)|=p$. If   $|Z(G)|=p^2$, then   $|Z(\mathrm{Inn}(G))|=p^2$ and therefore $|\gamma_2(G)/Z(G)|=p$. This implies that $|\gamma_2(G)|=p^3$ and $G/\gamma_2(G) \approx C_{p^2}\times C_{p}$.  Thus 
$|\mathrm{Aut}_z(G)|=|\mathrm{Hom}(C_{p^2}\times C_{p}, Z(G))|\geq p^3> |Z(\mathrm{Inn}(G))|,$ a contradiction to the assumption. This proves the `` {\em if} '' part of the theorem.

\section{\large Application}
In this section, we use the classification of all groups  of order $p^n$, where $p$ is an odd prime and $5\leq n\leq 6$, given by James [7].  As an application of our results, we find  those groups $G$ of order $p^5$ and $p^6$ for which $\mathrm{Aut}_z(G)=Z(\mathrm{Inn}(G))$. We shall mainly use the informations of these groups given in $\S 4.1$ and presentations of these groups  given in $\S 4.5$ and $\S 4.6$ of [7].

\begin{thm}
If $|G| = p^{5}$ and $cl(G)=3$, then 
$\mathrm{Aut}_z(G)=Z(\mathrm{Inn}(G))$
 if and only if $G$ is isomorphic to  $\Phi_{8}(32)$.
\end{thm}
{\bf Proof.}  There are only 2 isoclinism families {\em viz.} $\Phi_{7}$ and $\Phi_{8}$ which consist of groups $G$ such that $|Z(G)|=p$ and $cl(G)=3$. The family $\Phi_7$ consists of groups $G$ with $G/Z(G)\approx \Phi_2(1^4)$.   Thus  $G/Z(G)$ is of exponent $p$. Since $|\gamma_2(G)|=p^2$, it follows that $Z(G)< \gamma_2(G)$ and $G/\gamma_2(G)$ is an elementary abelian group of order $p^3$. Thus $d(G)=3$ and hence $\mathrm{Aut}_z(G)\ne Z(\mathrm{Inn}(G))$ by Theorem 3.2. The family $\Phi_8$ consists of only one group {\em viz.}
$$\Phi_{8}(32)=\;<\alpha_1,\alpha_2,\beta \;|\;[\alpha_1,\alpha_2]=\beta = \alpha_1^p,\;\beta^{p^2}=\alpha_2^{p^2}=1>$$ which is of rank 2 because $\beta = [\alpha_1,\alpha_2]\in \gamma_2(\Phi_{8}(32))\le \Phi(\Phi_{8}(32))$ and hence $\mathrm{Aut}_z(\Phi_{8}(32))=Z(\mathrm{Inn}(\Phi_{8}(32)))$
by Theorem 3.2.

\begin{thm}
If $|G|=p^{6}$ and $cl(G)=3$ or $4$, then
$\mathrm{Aut}_z(G)=Z(\mathrm{Inn}(G))$ if and only if 
$G$ is isomorphic to one of the groups in the isoclinism families
$\Phi_{25},\Phi_{26},\Phi_{28},\Phi_{29}$ and
$\Phi_{40}-\Phi_{43}$.
\end{thm}
{\bf Proof.} There are 16 isoclinism families {\em viz.} 
$\Phi_{22},\Phi_{24} - \Phi_{34}$ and $\Phi_{40} - \Phi_{43}$ which consist of groups $G$ for which $|Z(G)|=p$ and $cl(G)=3$ or 4. If $G$ is any group from the families $\Phi_{22},\Phi_{24},\Phi_{27}$ and $\Phi_{30} - \Phi_{33}$, then  $Z(G)< \gamma_2(G)$, $|\gamma_2(G)|=p^2$ or $p^3$ and $G/Z(G)$ is of exponent $p$. Therefore $G$ is of rank 3 or 4 and hence $\mathrm{Aut}_z(G)\ne Z(\mathrm{Inn}(G))$ by Theorem 3.3. The family $\Phi_{34}$ consists of groups $G$ such that $|\gamma_2(G)|=p^3$. From the presentations of these groups, it follows that they  are generated by 6 elements $\alpha,\alpha_1,\alpha_2,\beta_1,\beta_2,\gamma$ such that $\beta_1,\beta_2,\gamma,\alpha^p,\alpha_1^p,\alpha_2^p\in \gamma_2(G)$. This implies that  $\mho_1(G)\leq \gamma_2(G)$. Thus $d(G)=3$ and hence $\mathrm{Aut}_z(G)\ne Z(\mathrm{Inn}(G))$ by Theorem 3.3. Any group $G$ of the families $\Phi_{25},\;\Phi_{26},\;\Phi_{28}$ and $\Phi_{29}$ is generated by 5 elements $\alpha,\alpha_1,\alpha_2,\alpha_3,\alpha_4$ such that $\alpha_2,\alpha_3,\alpha_4\in \gamma_2(G)$. Thus $d(G)=2$ and hence $\mathrm{Aut}_z(G)= Z(\mathrm{Inn}(G))$ by Theorem 3.3.
If $G$ is any group from the families $\Phi_{40} - \Phi_{43}$, then $|\gamma_2(G)|=p^4$ and hence $d(G)=2$. Thus
$\mathrm{Aut}_z(G)=Z(\mathrm{Inn}(G))$ by Theorem 3.3.\\

\noindent {\bf Acknowledgement:} Authors are very thankful to the referee for his/her useful comments and suggestions.  

\begin{center} {\bf References}
\end{center}

\noindent [1] J.E. Adney and T. Yen, Automorphisms of a p-group, {\em Illinois J. Math.} {\bf 9} (1965), 137-43.

\noindent [2] Y. Berkovich, {\em Groups of Prime Power Order Volume 1}, (Walter de Gruyter, Berlin, 2008).

\noindent [3] N. Blackburn, On a special class of $p$-groups, {\em Acta Math.}, {\bf 100} (1958), 45-92.

\noindent [4] Y. Cheng, On double centralizer subgroups of some finite $p$-groups, {\em Proc. Amer. Math. Soc.}, {\bf 86} (1982), 205-208.

\noindent [5] M. J. Curran, Finite  groups   with  central automorphism   group  of  minimal  order, {\em Math. Proc. Roy. Irish Acad.}, {\bf 104A} (2) (2004),  223-229.

\noindent [6] M.J. Curran and D.J. McCaughan, Central  automorphisms   that   are   almost inner, {\em Comm. Algebra}, {\bf 29} (5) (2001), 2081-2087.

\noindent [7] R. James, The groups of order $p^{6}$ ($p$ an odd prime), {\em Math. Comp.}, {\bf 34} (150) (1980), 613-637.

\noindent [8] J. Wiegold, Multiplicators and groups with finite central factor-groups, {\em Math. Z.}, {\bf 89} (1965), 345-347.

\end{document}